\documentstyle[amstex, amscd,
12pt
]{article}

\def\ga{{\ifmmode \alpha \else $\alpha$ \fi }}
\def\gb{{\ifmmode \beta \else $\beta $ \fi }}
\def\gga{{\ifmmode \gamma \else $\gamma $ \fi }}
\def\gd{{\ifmmode \delta \else $\delta $ \fi }}
\def\gep{{\ifmmode \epsilon \else $\epsilon $ \fi }}
\def\gev{{\ifmmode \varepsilon \else $\varepsilon $ \fi }}
\def\gz{{\ifmmode \zeta \else $\zeta $ \fi }}
\def\get{{\ifmmode \eta \else $\eta $ \fi }}
\def\gth{{\ifmmode \theta \else $\theta $ \fi }}
\def\gtv{{\ifmmode \vartheta \else $\vartheta $ \fi }}
\def\gi{{\ifmmode \iota \else $\iota $ \fi }}
\def\gk{{\ifmmode \kappa \else $\kappa $ \fi }}
\def\gl{{\ifmmode \lambda \else $\lambda $ \fi }}
\def\gm{{\ifmmode \mu \else $\mu $ \fi }}
\def\gn{{\ifmmode \nu \else $\nu $ \fi }}
\def\gx{{\ifmmode \xi \else $\xi $ \fi }}
\def\gom{{\ifmmode \omicron \else $\omicron $ \fi }}
\def\gp{{\ifmmode \pi \else $\pi $ \fi }}
\def\gr{{\ifmmode \rho \else $\rho $ \fi }}
\def\gs{{\ifmmode \sigma \else $\sigma $ \fi }}
\def\gsv{{\ifmmode \varsigma \else $\varsigma $ \fi }}
\def\gt{{\ifmmode \tau \else $\tau $ \fi }}
\def\gu{{\ifmmode \upsilon \else $\upsilon $ \fi }}
\def\gph{{\ifmmode \phi \else $\phi $ \fi }}
\def\gpv{{\ifmmode \varphi \else $\varphi $ \fi }}
\def\gch{{\ifmmode \chi \else $\chi $ \fi }}
\def\gps{{\ifmmode \psi \else $\psi $ \fi }}
\def\go{{\ifmmode \omega \else $\omega $ \fi }}
\def\gG{{\ifmmode \Gamma \else $\Gamma $ \fi }}
\def\gD{{\ifmmode \Delta \else $\Delta $ \fi }}
\def\gTh{{\ifmmode \Theta \else $\Theta $ \fi }}
\def\gL{{\ifmmode \Lambda \else $\Lambda $ \fi }}
\def\gX{{\ifmmode \Xi \else $\Xi $ \fi }}
\def\gP{{\ifmmode \Pi \else $\Pi $ \fi }}
\def\gS{{\ifmmode \Sigma \else $\Sigma $ \fi }}
\def\gU{{\ifmmode \Upsilon \else $\Upsilon $ \fi }}
\def\gPh{{\ifmmode \Phi \else $\Phi $ \fi }}
\def\gPs{{\ifmmode \Psi \else $\Psi $ \fi }}
\def\gO{{\ifmmode \Omega \else $\Omega $ \fi }}

\def\BbA{{\ifmmode {\Bbb A}\else $\Bbb A$ \fi }}
\def\BbC{{\ifmmode {\Bbb C}\else $\Bbb C$ \fi }}
\def\BbF{{\ifmmode {\Bbb F}\else $\Bbb F$ \fi }}
\def\BbH{{\ifmmode {\Bbb H}\else $\Bbb H$ \fi }}
\def\BbN{{\ifmmode {\Bbb N}\else $\Bbb N$ \fi }}
\def\BbP{{\ifmmode {\Bbb P}\else $\Bbb P$ \fi }}
\def\BbQ{{\ifmmode {\Bbb Q}\else $\Bbb Q$ \fi }}
\def\BbR{{\ifmmode {\Bbb R}\else $\Bbb R$ \fi }}
\def\BbW{{\ifmmode {\Bbb W}\else $\Bbb W$ \fi }}
\def\BbZ{{\ifmmode {\Bbb Z}\else $\Bbb Z$ \fi }}

\def\2ff{{\ifmmode {{\text{\normalshape I}}@!@!@!@!{\text{\normalshape I}}} 
\else ${\text{\normalshape I}}@!@!@!@!{\text{\normalshape I}}$ \fi }}

\def\Boxtimes{{\ifmmode {\Box\!\!\!\!@!@!@!@!@!\times} 
\else $\Box\!\!\!\!@!@!@!@!@!\times$ \fi }}

\def\boxtimes{{\ifmmode ${\fbox{$\times$}}$ 
\else \fbox{$\times$} \fi }}

\newtheorem{thm}{Theorem}[section]
\newtheorem{lm}[thm]{Lemma}
\newtheorem{cor}[thm]{Corollary}
\newtheorem{prop}[thm]{Proposition}
\newtheorem{dfn}[thm]{Definition}
\newtheorem{con}[thm]{Conjecture}
\newtheorem{ntn}[thm]{Notation}
\newtheorem{propdef}[thm]{Proposition--Definition}
\newtheorem{remark}[thm]{Remark}

\newcommand{\bt}{\begin{thm}}
\newcommand{\et}{\end{thm}}
\newcommand{\bl}{\begin{lm}}
\newcommand{\el}{\end{lm}}
\newcommand{\bc}{\begin{cor}}
\newcommand{\ec}{\end{cor}}
\newcommand{\bp}{\begin{prop}}
\newcommand{\ep}{\end{prop}}
\newcommand{\bd}{\begin{dfn}}
\newcommand{\ed}{\end{dfn}}
\newcommand{\bcon}{\begin{con}}
\newcommand{\econ}{\end{con}}
\newcommand{\bntn}{\begin{ntn}}
\newcommand{\entn}{\end{ntn}}
\newcommand{\bpd}{\begin{propdef}}
\newcommand{\epd}{\end{propdef}}
\newcommand{\brem}{\begin{remark}}
\newcommand{\erem}{\end{remark}}

\newcommand{\beqn}{\begin{equation}}
\newcommand{\eeqn}{\end{equation}}
\newcommand{\pf}{\noindent {\bf Proof:}\hspace{2.5mm}}
\newcommand{\qed}{~\hspace{2.5mm}~\rule{2.5mm}{2.5mm}\vspace{\baselineskip}\\}

\newcommand{\hs}{\hspace{1.25mm}}
\newcommand{\hsp}{\hspace{2.5mm}}
\newcommand{\Hsp}{\hspace{5mm}}

\newcommand{\nim}{\text{\normalshape im}\hspace{.5mm}}
\newcommand{\nker}{\text{\normalshape ker}\hspace{.5mm}}

\newcommand{\ndim}{\text{\normalshape dim}\hspace{.5mm}}

\newcommand{\sym}{\text{\normalshape Sym}}
\newcommand{\Sym}{\text{\normalshape Sym}}

\newcommand{\Hom}{\text{\normalshape Hom}}

\begin{document}

\title{Hodge--Gaussian maps
\footnote{
\newline
The authors are members of CNR--GNSAGA (Italy). 
\newline
During the preparation of this paper the authors were partially supported by
National Research Project ``Geometria algebrica, algebra 
commutativa e aspetti computazionali'' of MURST--Italy. }}
\author{Elisabetta Colombo
\\Dipartimento di Matematica, Universit\`{a} di Milano
\and 
Gian Pietro Pirola
\\Dipartimento di Matematica, Universit\`{a} di Pavia
\and 
Alfonso Tortora
\\Dipartimento di Matematica, Universit\`{a} di Milano}
\date{
}
\maketitle


\begin{abstract}
\noindent Let $X$ be a compact K\"{a}hler manifold, and let $L$ be a line bundle
on $X.$ Define $I_k(L)$ to be the kernel of the multiplication map 
$\sym^kH^0(L)\rightarrow H^0(L^k).$ For all $h\leq k,$ we define a map
$$\rho:I_k(L)\rightarrow\text{Hom}(H^{p,q}(L^{-h}), H^{p+1,q-1}(L^{k-h})).$$
When $L=K_X$ is the canonical bundle, the map \gr computes a second fundamental 
form associated to the deformations of $X.$\\
If $X=C$ is a curve, then \gr is a lifting of the Wahl map
$I_2(L)\rightarrow H^0(L^2\otimes K_C^2).$\\
We also show how to generalize the construction of \gr to the cases of harmonic 
bundles and of couples of vector bundles.
\end{abstract}

\section*{Introduction}
In connection with the variations of Hodge structures (VHS), 
a number of authors have tackled the higher differentials of the period map.\\
A first definition of second fundamental form (2ff) for a  VHS of odd weight
is given in \cite{ivhs}. More recently, Karpishpan (\cite{k}) has defined a 2ff 
for VHS and
showed a way  to compute it for VHS coming from geometry, using Archimedean 
cohomology. In the case of curves, he asks whether this 2ff, at any given
point, lift $I_2(K_X)\rightarrow H^0(K_X^4),$ the second Wahl 
(or Gaussian) map for the canonical bundle.\\
In the projective case, second (and higher) fundamental
forms are defined for algebraic varieties, with respect to a fixed
projective embedding (cf. \cite{gh}, \cite{l}.)\\
In \cite{cime}, with reference to unpublished work of 
Green--Griffiths, it is reported that the  projective 2ff (in the sense of 
\cite{gh}) of any local Pl\"ucker embedding of the moduli space of curves gives,
as a quotient, the second Wahl map of the canonical 
line bundle.\\
Both kinds of 2ff, for VHS and projective embeddings, can be interpreted as 
instances of the (classical) 2ff $\2ff:=\pi\nabla|_{\cal S}$ associated to an 
extension of sheaves
$0\rightarrow{\cal S}\rightarrow E@>\pi>>{\cal Q}\rightarrow0,$
with $E$ a vector bundle with connection $\nabla.$\\
In this paper we define a family of maps, that we propose to call 
{\em Hodge--Gaussian maps,} existing  under very general 
conditions, namely for line bundles  over compact K\"{a}hler manifolds.
When applied to the canonical bundle,  the Hodge--Gaussian map is a  
2ff naturally associated to a  deformation of the manifold 
(see theorem \ref{hg2ff}.) If we are dealing with curves, we answer in the 
affirmative to the question asked in 
\cite{k}, consistently with the statement of \cite{cime}, cited earlier.
Actually, our result holds in a more general setup than those of both 
\cite{cime} and \cite{k}, in that it concerns not only the canonical bundle, 
but any line bundle on a curve. 
Also, the possibility of making explicit computations, at least in the case of 
curves, as in lemma \ref{lemmino}, seems to the authors a step 
towards understanding the curvature of the moduli space of curves.\\
The starting point for this paper was a construction of one of the authors 
(cf.~\cite{p}), that turned out to be a special case of ours. The hunch that it should be a kind of 2ff, and an attempt at understanding it as a lifting of a
Wahl map, in the spirit of Green--Griffiths, lead us to the present results.\\ 
The main idea underlying all of our maps is the following:\\
Let $L$ be a line bundle over a compact K\"{a}hler manifold $X,$ 
with $h^0(L)>1.$
Set  $I_2(L):=\ker(\sym^2H^0(L)\rightarrow H^0(L^2)).$ If 
$\gx=[\gth]\in H^1(L^{-1}),\hs \gth$ a Dolbeault representative of $\gx,$ 
and $\gl_i,\hs i=1,\dots,r,$ is a basis of $H^0(L),$ 
then the cup products $\gth\gl_i\in{\cal A}^{0,1}(X)$ 
have harmonic decompositions $\gth\gl_i=\gga_i+\bar{\partial}h_i.$ 
Now, for any $Q=\sum a_{ij}\gl_i\otimes\gl_j\in I_2(L),$ the section 
$$\sum a_{ij}\gl_i\partial h_j\in A^{1,0}(L)$$
determines an element of $H^0(L\otimes\gO^1_X).$\\
It turns out that the map 
$I_2(L)\otimes H^1(L^{-1})\rightarrow H^0(L\otimes\gO_X)$ is well
defined. Especially, when $X$ is a curve, this map,
seen as a map $I_2(L)\rightarrow H^0(L\otimes K_X)\otimes H^0(L\otimes K_X),$ 
is a lifting of the second Wahl map for 
$L, \hsp \gm_{_2}:I_2(L)\rightarrow H^0(L^2\otimes K_X^2),$
with respect to the natural multiplication map
$H^0(L\otimes K_X)\otimes H^0(L\otimes K_X)\rightarrow H^0(L^2\otimes K_X^2).$\\
The crux of our construction is the harmonic decomposition of the 
(p,q)--forms, to define the map, and the principle of two types, to prove
that it is well--defined.\\
This observation allows us to generalize the construction to a map 
$$I_k(L)\otimes H^{p,q}(L^{-h})\rightarrow H^{p+1,q-1}(L^{k-h}),$$
defined for line bundles $L$ on $X.$\\
Actually, the basic trick in the definition of the map is a switch from
$\bar{\partial}$ to $\partial,$ and it works also in 
more general situations, provided some kind of harmonic 
decomposition exist, for which the principle of two types holds.
This is the case for harmonic bundles, which admit the same kind of maps.
Such a generalization is not gratuitous, but with an eye 
towards finding interactions between Hodge theory and the equations defining 
an algebraic variety.\\
The authors' opinion is that the main interest of the present paper resides in 
the construction of a natural map \gr, not hiterto known in the literature. 
Indeed, in the published account \cite{cime} of the work of Green--Griffiths 
cited above, there is no mention of it.\\
Several people, whose encouragement we gratefully acknowledge, held the 
opinion that the non--holomorphic map \gr could be a suitable projection of an 
algebraic one. Its being non--holomorphic is likely to be the main obstruction 
to a more systematic use of \gr in algebraic geometry. However, a most likely 
application of \gr should be found in the investigation of the
curvature properties of certain moduli spaces, which fact nicely ties in with 
the non--holomorphicity. On the other hand, also the 2ff defined in \cite{k} 
is non--holomorphic even though it is somewhat shrouded in the use of 
Archimedean cohomology.\\ 
The paper is organized as follows:\\
In section 1 we define the Hodge--Gaussian map
$$\rho:I_k(L)\otimes H^{p,q}(L^{-h})\rightarrow H^{p+1,q-1}(L^{k-h})$$ 
whose construction is outlined above. We also note some formal properties of
the map, which are summarized in proposition \ref{sonomario}.\\
In section 2 we compare our map and the 2ff. Given a  smooth deformation  ${\cal X}@>\psi>>B$ of $X=X_{b_0},$ let $K_{{\cal X}|B}$ be the relative
canonical bundle. We show that the 2ff associated to the map
$\sym^k\psi_*K_{{\cal X}|B}\rightarrow\psi_*K_{{\cal X}|B}^k,$ at the point 
$b_{_0}\in B,$  is factorized by 
$I_k(K_X)\otimes H^{n-1,1}(K_X^{-1})@>\gr>>H^{n,0}(K_X^{k-1}),$ through the 
Kodaira--Spencer map 
$\gk:T_{B,b_0}\rightarrow H^1(T_X)\simeq H^{n-1,1}(K_X^{-1}).$\\
Section 3 deals with the case when $X=C$ is a curve:  we show 
that $\rho$ gives a lifting of the Wahl map.\\ 
In section 4 we show how to carry the construction of \gr over to more general 
situations, defining a Hodge--Gaussian map also in the following cases:
(a) for couples of vector bundles $E, F$--with $I_2(L)$ replaced by the second 
module of relations $R_2(E, F)$--and 
(b) for harmonic bundles.
\vspace\baselineskip\\
{\bf Acknowledgments:} The authors thank Eduard Looijenga, Marco Manetti, 
Eckart Viehweg and Claire Voisin for fruitful discussions on the topics of the 
present paper.

\section{The main construction}
Let $X$ be a compact K\"{a}hler manifold, $\ndim X=n,$ and let $L$ be a 
line bundle over $X, \hs h^0(L)=r>0.$\\ 
The goal of this section is to define the {\em Hodge--Gaussian map}
$$\gr:I_k(L)\rightarrow\text{Hom}(H^{p,q}(L^{-m}), H^{p+1,q-1}(L^{k-m})),$$
where 
$I_k(L):=\ker({\mathbf m}_k:\Sym^kH^0(L)\rightarrow H^0(L^k)),\hs {\mathbf m}_k$ 
being the multiplication map.\\
To do so, we need the following classical results of Hodge theory 
(see e.g. \cite{pag}, p. 84 and 149; also,
for a thorough exploitation of the principle of two types, \cite{dgms}.)
\bt\label{hodge}
Let $X$ be a compact K\"{a}hler manifold.\\
1. (Hodge theorem) Any $\bar{\partial}$-closed form $\ga\in A^{p,q}(X)$ has a 
unique harmonic representative, hence can be written as 
$\ga=\gga+\bar{\partial}h,$
with $\gga\in{\cal H}^{p,q}$  harmonic and $h\in A^{p,q-1}(X).$\\
2. (Principle of two types) Let $\ga\in A^{p,q}(X)$ satisfy 
$\partial\ga=\bar{\partial}\ga=0$ and be either $\partial-$ or 
$\bar{\partial}-$exact. Then for some $\gb\in A^{p-1,q-1}(X),$
$\ga=\partial\bar{\partial}\gb.$
\et
We now introduce some multi--index notation.\\
Fix a basis $\gl_{_1},\ldots,\gl_r$ of $H^0(L).$\\
Define $R_k:=\{1,2,\dots, r\}^k.$ If $S=(s_1,\dots,s_h)\in R_h$ and 
$T=(t_1,\dots,t_k)\in R_k,$ we denote 
$ST:=(s_1,\dots,s_h,t_1,\dots,t_k)\in R_{h+k}.$\\
\vspace{-1.5\baselineskip}
\begin{tabbing}
For any $J\in R_k$ we write: \= 
$a_{_J}=a_{j_1\dots j_k}\in\BbC$ is a scalar,\\
\> $\gl_{_{\otimes J}}=
\gl_{j_1}\otimes\dots\otimes\gl_{j_k}\in\otimes^kH^0(L),$\\
\>$\gl_{_J}=\gl_{j_1}\cdots\gl_{j_k}\in H^0(L^k).$\\
\end{tabbing}
\vspace{-1.5\baselineskip}
Clearly, an element $P\in I_k(L)$ is uniquely written as   
$\sum_{J\in R_k}a_{_J}\gl_{_{\otimes J}},$ with the $a_{_J}$'s symmetric in the $j$'s,  satisfying 
$\sum_{J\in R_k}a_{_J}\gl_{_J}=0.$\\
In standard multi--index notation, $P\in I_k(L)$ can be thought of as a 
polynomial of degree $k,\hsp \sum_{|K|=k}a_{_K}x^K,$ vanishing in \gl, i.e. 
$P(\gl)=\sum_{|K|=k}a_{_K}\gl^K=0,$ 
where $\gl^K=\gl_{_1}^{k_1}\cdots\gl_r^{k_r}.$\\  
\bpd\label{uno}
{\text{\bf : Hodge--Gaussian maps}}\\
Given $\gx\in H^{p,q}(L^{-m}),$ choose a Dolbeault representative  
$\gth\in A^{p,q}(L^{-m}).$ For any $T\in R_m,$ the cup product 
$\gth\gl_{_T}\in A^{p,q}(X)$ is $\bar{\partial}$-closed, so it has a harmonic 
decomposition
\begin{equation}
\gth\gl_{_T}=\gga_{_T}+\bar{\partial}h_{_T}\label{decomposition}
\end{equation}
with $\gga_{_T}\in{\cal H}^{p,q}$ and $h_{_T}\in A^{p,q-1}(X).$
Let $P=\sum_{J\in R_k}a_{_J}\gl_{_{\otimes J}}\in I_k(L).$\\ 
For all $0 < m\leq k,$ the following map is well-defined and $\BbC$--linear
$$
\begin{array}{cccc}
\gr: & I_k(L) & \rightarrow & 
\text{\normalshape Hom}(H^{p,q}(L^{-m}), H^{p+1,q-1}(L^{k-m}))\\
& P & \rightarrow & (\gx\mapsto\gr_{_P}(\gx))
\end{array}
$$
where $\gr_{_P}(\gx)$ is the Dolbeault cohomology class of the 
$L^{k-m}$--valued $(p+1,q-1)$-form
\begin{equation}
\gs_{_P}(\gth):=
\sum_{\stackrel{S\in R_{k-m}}{\scriptscriptstyle T\in R_m}}
a_{_{ST}}\gl_{_S}\partial h_{_T}.\label{last}
\end{equation}
\epd
\pf
 We need to check that $\gs_{_P}(\gth)$ is $\bar{\partial}$--closed and that 
$\gr_{_P}(\gx)=[\gs_{_P}(\gth)],$ as an element of $H^{p+1,q-1}(L^{k-m})$, 
is independent of the choices made.\\
(i)\Hsp $\gs_{_P}(\gth)$ is $\bar{\partial}$--closed.\\
Indeed, 
\begin{eqnarray}
\bar{\partial}\gs_{_P}(\gth)&=& 
\sum a_{_{ST}}\gl_{_S}\bar{\partial}\partial h_{_T}=
-\sum a_{_{ST}}\gl_{_S}\partial(\gth\gl_{_T}-\gga_{_T})\nonumber\\
&=&-\sum a_{_{ST}}\gl_{_S}\partial(\gth\gl_{_T})\nonumber
\end{eqnarray}
because $\gga_{_T}$ is harmonic. A local computation 
shows that $\sum a_{_{ST}}\gl_{_S}\partial(\gth\gl_{_T})$ vanishes: for any 
$p\in X,$ let $\ell$ and $\ell^*$ be a local generator of $L$ and its dual in a 
neighborhood $U$ of $p,$ then 
$\gl_{_i}=\phi_{_i}\cdot\ell,\;\gth=\tau\cdot(\ell^*)^m,$ 
with $\phi_{_i}$ functions and $\tau$ a (p,q)--form on $U$ respectively. 
On $U,$ we have $\gl_{_T}=\phi_{_T}\ell^m,$ where 
$\phi_{_T}=\phi_{_{t_1}}\cdot\ldots\cdot\phi_{_{t_m}}$
is a function defined on $U,$ so $\partial(\gth\gl_{_T})=
\phi_{_T}\partial\tau+(-1)^{p+q}\tau\wedge\partial\phi_{_T},$ hence 
$$\sum a_{_{ST}}\gl_{_S}\partial(\gth\gl_{_T})=
\left(\partial\tau\sum a_{_{ST}}\phi_{_S}\phi_{_T}+(-1)^{p+q}\tau\wedge
\sum a_{_{ST}}\phi_{_S}\partial\phi_{_T}\right)\ell^{k-m}=0.$$
Indeed, $\sum a_{_{ST}}\gl_{_S}\gl_{_T}=0$ means that the function 
$\sum a_{_{ST}}\phi_{_S}\phi_{_T}$ is identically zero on $U,$ thus also
$\partial\left(\sum a_{_{ST}}\phi_{_S}\phi_{_T}\right)=0;$ since the scalars
$a_{_J}$ are symmetric with respect to the indices $j$'s, it is easy to see 
that $\sum a_{_{ST}}\phi_{_S}\partial\phi_{_T}=
\frac{m}{k}\partial(\sum a_{_{ST}}\phi_{_S}\phi_{_T})=0.$ 
(see {\it infra}, remark \ref{green})\\
(ii)\Hsp $\gr_{_P}(\gx)$ does not depend on the choice of $\gth.$\\ 
Let $\tilde{\gth}$ be another Dolbeault representative of \gx, we have
$\tilde{\gth}=\gth+\bar{\partial}\chi,$ with $\chi\in A^{p,q-1}(L^{-m}).$ Now, 
$$\tilde{\gth}\gl_{_T}=\gga_{_T}+\bar{\partial}\tilde{h}_{_T},$$ 
with $\gga_{_T}$ unchanged
because it is the unique harmonic representative of
$\gx\gl_{_T}\in H^{p,q}(X),$ so 
$$\gga_{_T}+\bar{\partial}\tilde{h}_{_T}=
\tilde{\gth}\gl_{_T}=(\gth+\bar{\partial}\chi)\gl_{_T}=\gth\gl_{_T}+
\bar{\partial}(\chi\gl_{_T})=
\gga_{_T}+\bar{\partial}h_{_T}+\bar{\partial}(\chi\gl_{_T}),$$ 
hence
$$\bar{\partial}\tilde{h}_{_T}=\bar{\partial}(h_{_T}+\chi\gl_{_T})$$ 
and
$$\tilde{h}_{_T}=h_{_T}+\chi\gl_{_T}+g_{_T},$$
with $g_{_T}$ a $\bar{\partial}$--closed $(p,q-1)$--form.\\
It follows 
$$\sum a_{_{ST}}\gl_{_S}\partial\tilde{h}_{_T}=
\sum a_{_{ST}}\gl_{_S}\partial h_{_T}+\sum a_{_{ST}}\gl_{_S}\partial g_{_T},$$
because $\sum a_{_{ST}}\gl_{_S}\partial (\chi\gl_{_T})=0$ as above. 
So we need to show that $\sum a_{_{ST}}\gl_{_S}\partial g_{_T}$ 
is $\bar{\partial}$--exact.\\
As $\partial g_{_T}$ is $\bar{\partial}$--closed, by the principle of two types 
$\partial g_{_T}=\bar{\partial}\partial k_{_T},$ hence
$$\sum a_{_{ST}}\gl_{_S}\partial g_{_T}=\sum a_{_{ST}}\gl_{_S}\bar{\partial}\partial k_{_T}=
\bar{\partial}\left(\sum a_{_{ST}}\gl_{_S}\partial k_{_T}\right).$$
The linearity of $\gr_{_P}$ and its independence of the choice of a basis of
$H^0(L)$ are clear.\qed
\brem\label{green}
{\em (i) $\gr_{_P}(\gx)$ can also be defined, perhaps more 
intuitively, thinking of $P\in I_k(L)$  as a polynomial of degree $k$ 
vanishing in \gl, $P(\gl)=0.$
If we write $P=P(x)$ in the form $\sum_{J\in R_k}a_{_J}x_{_J},$ where $a_{_J}$'s
are the same scalar seen above and $x_{_J}=x_{_{j_1}}\dots x_{_{j_k}},$ then 
it is easy to see that the partial derivatives of $P(x)$ are given by 
$\frac{\partial P}{\partial x_i}=k\sum_{S\in R_{k-1}}a_{_{Si}}x_{_S}.$
Also, when $\gx\in H^{p,q}(L^{-1}),$ (\ref{decomposition}) becomes 
$\gth\gl_i=\gga_{_i}+\bar{\partial}h_{_i},$ for all $i=1,\dots,r.$ 
Thus $\gr_{_P}(\gx)$ is the cohomology class of the form
$$\gs_{_P}(\gth)=
\frac{1}{k}\sum_{i=1}^{r}\frac{\partial P}{\partial x_i}(\gl)\partial h_{_i}.$$
For $m>1,$ the formula expressing $\gr_{_P}(\gx)$ in terms of higher--order 
derivatives of $P(x)$ is slightly more complicated. For all 
$T=(t_1,\dots,t_m)\in R_m,$ let
$\partial_{_T}P=\frac{\partial^m P}{\partial x_{t_1}\dots\partial x_{t_m}},$
then one sees that
$\partial_{_T}P=\frac{k!}{(k-m)!}\sum_{S\in R_{k-m}}a_{_{ST}}x_{_S},$
so $\gs_{_P}(\gth)=\sum_{S,T}a_{_{ST}}\gl_{_S}\partial h_{_T}=
\frac{(k-m)!}{k!}\sum_T\partial_{_T}P(\gl)\partial h_{_T}.$ Now, in standard 
multiindex notation, $\partial_{_T}P=\frac{\partial^m P}{\partial x^{\gt}},$ 
where $\gt=\gt(T)=(\gt_{_1},\ldots,\gt_{_r}),\hs \gt_{_j}$ being how many times 
$j$ appears in $T=(t_1,\dots,t_m).$ Also, the same derivative 
$\frac{\partial^m P}{\partial x^{\gt}}$ is repeated 
$\frac{m!}{\gt_1!\cdots\gt_r!}$ times, corresponding to the different 
$T\in R_m$ which give the same $\gt(T).$ Summing up, we obtain 
$$\gs_{_P}(\gth)=\frac{m!(k-m)!}{k!}\sum_{|I|=m}\frac{1}{i_1!\cdots i_r!}
\frac{\partial^mP}{\partial x^I}(\gl)\partial h^I,$$
with $h^I$ given by the decomposition (\ref{decomposition}) relative to
$\gth\gl^I=\gth\gl_{1}^{i_1}\cdots\gl_{r}^{i_r}.$\\
(ii) The basic trick in the definition of $\gr_{_P}$ is to take the  
$\bar{\partial}$--exact part of the decomposition of a form and then switch 
to a $\partial$--exact form, i.e. going from $\bar{\partial}h_{_T}$ to 
$\partial h_{_T}.$ To do so, we just need the two facts of 
theorem~\ref{hodge}, hence a similar construction can be carried out also in 
other more general situations, where we have some kind of harmonic 
decomposition, for which the principle of two types holds.}
\erem
\bp\label{kaler} 
If $X$ is a compact complex manifold having several 
K\"{a}hler metrics compatible with its complex structure, 
then the map $\gr_{_P}$ is independent of the 
(K\"{a}hler) metric used to define it, and is completely determined by the 
underlying complex structure of $X.$
\ep
\pf Let $K_1$ and $K_2$ be the harmonic projectors coming from two different
K\"{a}hler metrics on $X;$ then, for any  $\bar{\partial}$--closed form \go we 
have the harmonic decompositions $\go=K_i\go+\bar{\partial}h_i,\hs i=1, 2.$
Set $\psi:=\partial(h_1-h_2)=\partial h.$\\
We claim that \gps is $\bar{\partial}$--exact.\\
Indeed, $\bar{\partial}h$ is $\partial$--closed, 
(because $\bar{\partial}h=K_2\go-K_1\go,$ with the $K_i\go$ harmonic forms) 
so the  principle of two types implies that $\bar{\partial}h=-\bar{\partial}\partial f,$ or, equivalently, 
$\bar{\partial}(h+\partial f)=0.$ Therefore, $h+\partial f$ 
has harmonic decomposition 
$h+\partial f=\gt+\bar{\partial}l.$\\
It follows that
$\gps=\partial h=\partial(h+\partial f)=\partial\gt+\partial\bar{\partial}l=   
\bar{\partial}(-\partial l).$\\
Going back to our situation, $\gth\gl_{_T}$ has harmonic decompositions, 
with respect to the different K\"{a}hler structures, 
$\gth\gl_{_T}=\gga_{_T}+\bar{\partial}h_{_T}=\gd_{_T}+\bar{\partial}g_{_T},$
hence $\partial(h_{_T}-g_{_T})=\bar{\partial}l_{_T}$ is $\bar{\partial}$--exact.
It follows that 
$$\sum a_{_{ST}}\gl_{_S}\partial h_{_T}-\sum a_{_{ST}}\gl_{_S}\partial g_{_T}=
\sum a_{_{ST}}\gl_{_S}\bar{\partial}l_{_T}=
\bar{\partial}\left(\sum a_{_{ST}}\gl_{_S}l_{_T}\right),$$
thus the cohomology classes
$$\left[\sum a_{_{ST}}\gl_{_S}\partial h_{_T}\right]=
\left[\sum a_{_{ST}}\gl_{_S}\partial g_{_T}\right]$$
are equal in $H^{p+1,q-1}(L^{k-m}).~\hspace{2.5mm}~\rule{2.5mm}{2.5mm}$
\brem
{\em $\gr_{_P}$ does not vary holomorphically on family of varieties, in 
the following sense.\\
Let ${\cal X}\rightarrow S$ be a smooth analytic family of K\"{a}hler manifolds 
and let ${\cal L}\rightarrow{\cal X}$ be a line bundle. Define
${\cal H}^{p,q}({\cal L}^{k}):=R^q\pi_*(\gO^p_{_{{\cal X}|S}}\otimes L^k)$ 
and 
${\cal I}_r({\cal L}):=
{\ker}(\Sym^r\pi_*{\cal L}\rightarrow\pi_*{\cal L}^r).$ 
\gr extends to a map 
$$\tilde{\gr}:{\cal I}_k({\cal L})\otimes{\cal H}^{p,q}({\cal L}^{-m})
\rightarrow{\cal H}^{p+1,q-1}({\cal L}^{k-m})$$
which is not holomorphic, but only real--analytic.}
\erem
The maps $\gr$ have a few more properties worth noting. 
\bp\label{lineare}
For all $\ga\in H^{s,t}(X)$ and $\gx\in H^{p,q}(L^{-m}),$
$$\gr_{_P}(\ga\cdot\gx)=\ga\cdot\gr_{_P}(\gx).$$
\ep
\pf Choose a harmonic representative \gb of the class \ga, then, recalling the
notation of (\ref{decomposition}),
$\gb\wedge\gga_{_T},$ which is not harmonic, has a harmonic decomposition
$$\gb\wedge\gga_{_T}=\gs_{_T}+\bar{\partial}g_{_T}.$$  
The forms $\gb, \gga_{_T}$ and $\gs_{_T},$ being harmonic, are both 
$\partial-$ and  $\bar{\partial}-$closed, hence $\bar{\partial}g_{_T}$ 
is $\partial-$closed, so, by the principle of two types, 
$\bar{\partial}g_{_T}=\bar{\partial}\partial k_{_T},$
for a suitable $k_{_T}.$ It follows that 
$\gb\wedge\gth\gl_{_T}=\gb\wedge(\gga_{_T}+\bar{\partial}h_{_T})=
\gb\wedge\gga_{_T}+(-1)^{s+t}\bar{\partial}(\gb\wedge h_{_T})=
\gs_{_T}+\bar{\partial}\partial k_{_T}
+(-1)^{s+t}\bar{\partial}(\gb\wedge h_{_T})=
\gs_{_T}+\bar{\partial}(\partial k_{_T}+(-1)^{s+t}\gb\wedge h_{_T}),$
so we have the harmonic decomposition
$$\gb\wedge\gth\gl_{_T}=\gs_{_T}+\bar{\partial}f_{_T},$$
with $f_{_T}=\partial k_{_T}+(-1)^{s+t}\gb\wedge h_{_T},$ hence 
$\partial f_{_T}=\gb\wedge\partial h_{_T}$---recall that \gb is
$\partial-$closed. The outcome is that 
\begin{eqnarray}
\gr_{_P}(\ga\cdot\gx)&=&\left[\sum a_{_{ST}}\gl_{_S}\partial f_{_T}\right]=
\left[\sum a_{_{ST}}\gl_{_S}\gb\wedge\partial h_{_T}\right]\nonumber\\
&=&\left[\gb\wedge\sum a_{_{ST}}\gl_{_S}\partial h_{_T}\right]=
[\gb]\cdot\left[\sum a_{_{ST}}\gl_{_S}\partial h_{_T}\right]=\nonumber\\
&=&\ga\cdot\gr_{_P}(\gx)~\hspace{2.5mm}~\rule{2.5mm}{2.5mm}\nonumber
\end{eqnarray}
\bp\label{symmetry}
For all $P\in I_k(L), \hs \gx\in H^{p,q}(L^{-m})$ and 
$\eta\in H^{n-p-1,n-q+1}(L^{-k+m}),$
$$\gx\gr_{_P}(\get)=(-1)^{p+q+1}\get\gr_{_P}(\gx).$$
\ep
\pf Let $\gth\in A^{p,q}(L^{-m})$ and $\chi\in A^{n-p-1,n-q+1}(L^{-k+m})$ be
Dolbeault representatives of \gx and \get respectively. Given 
$S\in R_{k-m}, T\in R_m,$ consider the corresponding harmonic decompositions 
$\gth\gl_{_T}=\gga_{_T}+\bar{\partial}h_{_T}$ and 
$\chi\gl_{_S}=\gd_{_S}+\bar{\partial}k_{_S}.$ Then the cohomology class 
$\get\gr_{_P}(\gx)\in H^{n,n}(X)$ has Dolbeault representative
$$\chi\cdot\sum a_{_{ST}}\gl_{_S}\partial h_{_T}=
\sum a_{_{ST}}\chi\gl_{_S}\wedge\partial h_{_T}=
\sum a_{_{ST}}\gd_{_S}\wedge\partial h_{_T}+
\sum a_{_{ST}}\bar{\partial}k_{_S}\wedge\partial h_{_T}.$$
It is easy to see that $d(\gd_{_S}\wedge h_{_T})=\gd_{_S}\wedge\partial h_{_T},$
hence $\get\gr_{_P}(\gx)$ is represented also by the form  
$\sum a_{_{ST}}\bar{\partial}k_{_S}\wedge\partial h_{_T}.$ Similarly, 
$\gx\gr_{_P}(\get)$ is represented by 
$\sum a_{_{TS}}\bar{\partial}h_{_T}\wedge\partial k_{_S}.$\\
Now, if $h\in A^{p,q-1}(X)$ and $k\in A^{n-p-1,n-q}(X),$ 
it is true in general that 
$[\bar{\partial}h\wedge\partial k]=
(-1)^{p+q+1}[\bar{\partial}k\wedge\partial h]$ 
in $H^{n,n}(X).$\\
Indeed, taking into account the number of $dz$'s and $d\bar{z}$'s, one sees that
$d(h\wedge dk)=\partial h\wedge\bar{\partial}k+\bar{\partial}h\wedge\partial k,$
thus $[\bar{\partial}h\wedge\partial k]=-[\partial h\wedge\bar{\partial}k]$
in $H^{n,n}(X).$ Since $\partial h\wedge\bar{\partial}k=
(-1)^{(p+q)(2n-p-q)}\bar{\partial}k\wedge\partial h=
(-1)^{p+q}\bar{\partial}k\wedge\partial h,$ then 
$[\bar{\partial}h\wedge\partial k]=
(-1)^{p+q+1}[\bar{\partial}k\wedge\partial h].$\\
The conclusion is now clear---recall that the $a_{_J}$'s are symmetric with 
respect to the indices $j$'s:\\
\begin{eqnarray}
\gx\gr_{_P}(\get)&=&
\left[\sum a_{_{TS}}\bar{\partial}h_{_T}\wedge\partial k_{_S}\right]=
(-1)^{p+q+1}
\left[\sum a_{_{ST}}\bar{\partial}k_{_S}\wedge\partial h_{_T}\right]\nonumber\\
&=&(-1)^{p+q+1}\get\gr_{_P}(\gx)~\hspace{2.5mm}~\rule{2.5mm}{2.5mm}\nonumber
\end{eqnarray}
\bntn
Given a line bundle $L$ as before, write
$$H^{\bullet}(L^{\bullet}):=\oplus_{p,q,k}H^{p,q}(L^k),$$
with $p, q, k\in\BbN_0,\; p, q\leq n.$
\entn
Furthermore, standard notations are
$$H^{\bullet}(X):=\oplus_{p,q}H^{p,q}(X)\;\; \text{\normalshape and}\;\;
I(L):=\oplus_kI_k(L).$$
Clearly, $H^{\bullet}(L^{\bullet})$ has a structure of $H^{\bullet}(X)-$module, 
given by the cup product. 
Using the identification $H^{p,q}(L^{-m})^*=H^{n-p,n-q}(L^{m}),$ the map 
$\gr_{_P}\in\\ \text{\normalshape Hom}(H^{p,q}(L^{-m}), H^{p+1,q-1}(L^{k-m}))$ 
is an element of $H^{n-p,n-q}(L^m)\otimes_\BbC H^{p+1,q-1}(L^{k-m})),$
hence $\gr$ is a map 
$$\gr:I_k(L)\rightarrow H^{n-p,n-q}(L^m)\otimes_\BbC H^{p+1,q-1}(L^{k-m})).$$
Note that, when $L$ is ample, \gr is nonzero only when $p+q=n.$\\
Putting the $\gr$'s together, for all values of $k,$ and taking into account 
the linearity expressed in proposition \ref{lineare}, we have a map 
$$\rho:I(L)\rightarrow 
H^{\bullet}(L^{\bullet})\otimes_{H^{\bullet}(X)}H^{\bullet}(L^{\bullet}).$$
Thinking of $I(L)$ and 
$H^{\bullet}(L^{\bullet})\otimes_{H^{\bullet}(X)}H^{\bullet}(L^{\bullet})$ 
as graded $\BbC$--modules, via
$$I(L)=\oplus I_k(L), \;\;
H^{\bullet}(L^{\bullet})\otimes_{H^{\bullet}(X)}H^{\bullet}(L^{\bullet})=
\oplus_k\left(\oplus_{m+j=k}
H^{\bullet}(L^m)\otimes_{H^{\bullet}(X)}H^{\bullet}(L^j)\right),$$
\gr is then a map of graded modules.\\
Proposition \ref{symmetry} expresses the fact that $\gr(P)\in 
H^{\bullet}(L^{\bullet})\otimes_{H^{\bullet}(X)}H^{\bullet}(L^{\bullet})$ 
is invariant with respect to the involution
$$\iota(\xi\otimes\eta):
=(-1)^{{\text deg}\xi\cdot{\text deg}\eta+1}\eta\otimes\xi,$$
where $\deg\xi=p+q$ for $\xi\in H^{p,q}(L^k).$\\
We can summarize the remarks above in the following
\bp\label{sonomario}
$\rho:I(L)\rightarrow 
H^{\bullet}(L^{\bullet})\otimes_{H^{\bullet}(X)}H^{\bullet}(L^{\bullet})$
is a map of graded $\BbC$--modules. Its image is contained in the subspace of
$H^{\bullet}(L^{\bullet})\otimes_{H^{\bullet}(X)}H^{\bullet}(L^{\bullet})$
invariant with respect to the involution $\iota(\xi\otimes\eta)
=(-1)^{{\text deg}\xi\cdot{\text deg}\eta+1}\eta\otimes\xi.
~\hspace{2.5mm}~\rule{2.5mm}{2.5mm}$
\ep
Especially, when $n=\ndim X$ is odd, $n=2m+1,$ the map $\rho$ is 
symmetric on the middle cohomology, i.e. 
$$\rho:I_{2k}(L)\rightarrow\sym^2 H^{m+1,m}(L^{k}).$$

\section{Hodge--Gaussian map and\\ 
second fundamental form}
Let $X$ be a complex manifold and let $E$ be a holomorphic vector bundle on $X,$
with connection 
$\nabla:{\cal A}^0(E)\rightarrow{\cal A}^1(E).$ For any  exact sequence of
sheaves of ${\cal O}_X$--modules
$0\rightarrow {\cal S}\rightarrow E@>\pi>>{\cal Q}\rightarrow0,$ 
the second fundamental form (2ff) of ${\cal S}$ in $E$ is the 
${\cal A}^0(X)$--linear map 
$$\2ff:{\cal A}^0({\cal S})\rightarrow{\cal A}^1({\cal Q})$$
defined by $\2ff(\gs):=\pi\nabla|_{_{\cal S}}(\gs).$ 
If $\nabla$ is compatible with the complex structure, then \2ff lands into 
${\cal A}^{1,0}({\cal Q}),$ hence
$\2ff\in A^{1,0}({\mathit Hom}({\cal S},{\cal Q}))$ (see e.g. \cite{pag}.)
With an eye on the case at hand, we slightly enlarge the definition of 2ff by
allowing an exact sequence of type 
$0\rightarrow {\cal S}\rightarrow E@>\pi>>{\cal Q},$ for which \gp is not
necessarily surjective; clearly, the same definition of $\2ff$ makes 
still sense.\\
We are most interested in the following situation:\\
Let ${\cal X}@>\psi>>B$ be a smooth analytic family of compact K\"{a}hler 
manifolds of dimension $n,$ i.e. assume that \gps is a submersion, 
so all fibers are smooth. For all $k,$ there are exact sequences
\begin{equation}
0\rightarrow{\cal I}_k(K_{{\cal X}|B})
\rightarrow\sym^k(\psi_*K_{{\cal X}|B})
@>{\mathbf m}>>\psi_*K_{{\cal X}|B}^k,\label{famiglia}
\end{equation}
whose 2ff we denote by $\2ff_k.$ Here $K_{{\cal X}|B}$ is the relative canonical 
bundle, $K_{{\cal X}|B}=\wedge^n\gO^1_{{\cal X}|B}, \hsp\gO^1_{{\cal X}|B}$
being the relative cotangent bundle. Also, ${\mathbf m}$ is the natural 
multiplication map.\\ 
Recall that the fiber of 
$\psi_*K_{{\cal X}|B}$ on the point $b\in B$ is $H^0(X_b,K_{X_b}),$ so 
$\psi_*K_{{\cal X}|B}={\cal F}^n\subseteq{\bf R}^n\psi_*\BbC$ 
is a piece of the Hodge filtration and, as such, has a natural metric connection
induced by $\nabla^{\mathit\text{GM}},$ the flat Gauss--Manin (GM) connection 
on the polarized VHS\hsp ${\bf R}^n\psi_*\BbC.$
The 2ff of the exact sequence (\ref{famiglia})
$$\2ff_k:{\cal I}_k\rightarrow\psi_*K_{{\cal X}|B}^k\otimes\gO_B,$$
becomes, on the central fiber, 
$$\2ff_k:I_k\otimes T_{B,b_0}\rightarrow H^0(K_X^k).$$
Now, the Kodaira--Spencer (KS) map of the family ${\cal X}$ is 
(on the central fiber)
$$\gk:T_{B,b_0}\rightarrow H^1(T_X).$$
With the identifications 
$H^0(K_X^k)=H^0(K_X^{k-1}\otimes\gO_X^n)=H^{n,0}(K_X^{k-1})$ and
$H^1(T_X)=H^1(\gO^*_X)=H^1(\gO^{n-1}_X\otimes K_X^{-1})=
H^{n-1,1}(K_X^{-1}),$ we have the following statement.
\bt\label{hg2ff} 
The diagram
\[
\begin{CD}
I_k\otimes T_{B,b_0}@>{id\otimes\gk}>>I_k\otimes H^{n-1,1}(K_X^{-1})\\
@V{\2ff_k}VV @VV{\gr}V\\
H^{n,0}(K_X^{k-1})@= H^{n,0}(K_X^{k-1}) 
\end{CD}
\]
is commutative up to a constant.
\et
The strategy of proof is very simple. First, it is enough to consider only the 
case of one--dimensional deformations ${\cal X}\rightarrow\gD,$ hence we need
to check the equality for just one vector ${\bf v}\in T_{\gD,0}.$
We now compute both the 2ff and the KS map using a fixed 
$C^\infty$--trivialization ${\cal X}\simeq\gD\times X.$ Finally, for any given
$P\in I_k,$ we plug in the value $\gk({\bf v})$ in the expression of $\gr_{_P}$
and get $\2ff_{k,{\bf v}}(P)=\gr_{_P}(\gk({\bf v})),$ up to a constant.\\
Of course, there is a slight abuse of notation in denoting with \gr the map of 
the present theorem, map which is actually \gr only after an obvious duality.
\vspace{\baselineskip}\\
\pf Since the statement is local, we can suppose that ${\cal X}$ be a 
one--dimensional deformation, i.e. ${\cal X}@>\psi>>\gD$ 
be parameterized by the unit circle $\gD=\{|t|<1\},$ with $X_0=X$ and 
${\bf v}=\frac{\partial}{\partial t}\in T_{\gD,0}.$
The Hodge bundle ${\cal H}^n_{_\BbC}={\bf R}^n\psi_*\BbC$ is a flat bundle,
with a flat connection $\nabla^{\mathit\text{GM}},$ 
the (local) GM connection.
The GM connection induces the connection $\nabla^{n,0}$ on the subbundle 
${\cal H}^{n,0}=F^n{\cal H}^n_{_\BbC}=\psi_*K_{{\cal X}|\gD},$ 
which in turn induces a connection, 
denoted by $\nabla,$ on the symmetric product $\sym^k{\cal H}^{n,0}.$ \\
We now compute both the GM connection and the KS map following the method set 
forth in \cite{cime} pp. 30--32, which we briefly summarize here.\\
Let $Y$ be a $C^\infty$--lifting of the holomorphic vector field
$\frac{\partial}{\partial t}$ on \gD; then we get a $C^\infty$--trivialization 
$\gD\times X@>\gt>>{\cal X}$ by $\gt(x,t):=\Phi_{tY}(1),$ where $\Phi_{Y}(t)$
denotes the flow associated to the vector field $Y.$ \\
One sees that $\bar{\partial}Y|_X\in A^{0,1}(T_{\cal X}|_X)$ is actually a 
closed form $\gth=\bar{\partial}Y|_X$ in $A^{0,1}(T_X)$ that represents the KS 
class associated to $\frac{\partial}{\partial t},$ i.e. 
$\gk(\frac{\partial}{\partial t})=[\bar{\partial}Y|_X]=[\gth].$\\
Let $\go(t)$ be a section of ${\cal H}^{n,0},$ then, for all 
$t, \go(t)\in H^0(K_{X_t});$ we may think of $\go(t)$ as 
$\gO\in A^{n,0}({\cal X})$ such that $\gO|_{X_t}=\go(t)$ as $(n,0)$--forms on 
$X_t.$ \\
The isomorphism $\gt_{_t}:X\rightarrow X_{_t},$ induced by \gt, 
gives an inclusion $\gt_{_t}^*: A^{n,0}(X_{_t})\hookrightarrow A^n(X).$ 
Since $\go(t)\in A^{n,0}(X_t)$ is $d$--closed, so is also 
$\gt^*_t(\go(t))\in A^n(X),$ thus we obtain a power series expansion  around 
$t=0$
$$\gt^*_t(\go(t))=\go+(\ga+dh)t+\circ(t^2),$$
with $\go=\go(0), \ga$ a harmonic $n$--form and  $h$ an $(n-1)$--form.\\
It follows that, as cohomology classes, 
$\nabla^{\mathit\text{GM}}_{\frac{\partial}{\partial t}}[\go(t)]_{_{t=0}}=
[\ga].$\\
On the other hand, $\frac{\partial\go(t)}{\partial t}=\gt^*{\cal L}_Y\gO=
\gt^*<d\gO,Y>+\gt^*d<\gO,Y>$ as forms, so $\frac{\partial\go(t)}{\partial t}$
has at least $n-1\, dz$'s, hence $\ga+dh$ lives in 
$A^{n,0}(X)\oplus A^{n-1,1}(X)$ and is of type
$\ga=\ga^{n,0}+\ga^{n-1,1}, \hsp h\in A^{n-1,0}(X).$\\
Finally, the $(n-1,1)$ part of $\frac{\partial\go(t)}{\partial t}|_{t=0}$ is the 
contration of \gO with $\bar{\partial}Y,$ restricted to $X.$ So we have the 
harmonic decomposition
\begin{equation}
\gth\go=\ga^{n-1,1}+\bar{\partial}h,\label{su}
\end{equation}
and $\gk(\frac{\partial}{\partial t})\cdot\go=[\ga^{n-1,1}].$\\
Now the conclusion of the proof is straightforward.\\
Let $P\in I_k,$ it is of type $P=\sum a_{_J}\go_{_{\otimes J}},$ so 
$\2ff_{k,{\frac{\partial}{\partial t}}}(P)=
{\mathbf m}({\nabla}_{\frac{\partial}{\partial t}}\gs|_{_{t=0}}),$ 
where $\gs(t)$ is a section through $P=\gs(0).$ Since 
${\nabla}^{n,0}_{\frac{\partial}{\partial t}}\go_{_j}(t)|_{_{t=0}}=
\ga_{_j}^{n,0},$ we see that 
$\2ff_{k,{\frac{\partial}{\partial t}}}(P)$ is represented by the form 
$\sum\dot{a}_{_J}(0)\go_{_ J}+
k\sum_{\stackrel{j=1\ldots g}{\scriptscriptstyle L\in R_{k-1}}}
a_{_{Lj}}(0)\go_{_L}\ga_{_j}^{n,0}$ ---here $g=h^0(K_X);$ recall that 
$R_{k-1}=\{ 1,\ldots,g\}^{k-1},$ also see remark \ref{green}.\\ 
Since $\gs(t)$ is a section of\, ${\cal I}_k,$ we have that
$\sum a_{_J}(t)\go(t)_{_J}=0$ identically, so also its derivative with 
respect to $t$ vanishes at $t=0,$ i.e. 
$\sum\dot{a}_{_J}(0)\go_{_ J}+k\sum a_{_{Lj}}(0)\go_{_L}(\ga_{_j}+dh_{_j})=0,$
and, taking the $(n,0)$ part of $\ga_{_j}+dh_{_j},\hsp
\sum\dot{a}_{_J}(0)\go_{_ J}+k\sum a_{_{Lj}}(0)\go_{_L}\ga_{_j}^{n,0}=
-k\sum a_{_{Lj}}(0)\go_{_L}\partial h_{_j}.$\\
In other words, up to a constant factor,
$$\2ff_{k,{\frac{\partial}{\partial t}}}(P)=
\sum a_{_{Lj}}\go_{_L}\partial h_{_j},$$
where $a_{_{Lj}}=a_{_{Lj}}(0).$\\
To compute $\gr_{_P}(\gk({\frac{\partial}{\partial t}}))$ we take \gth as 
representative of $\gk({\frac{\partial}{\partial t}}),$ so we have the harmonic 
decompositions (\ref{su}) relative to the products $\gth\go_{_j},$ i.e.
$\gth\go_{_j}=\ga_{_j}^{n-1,1}+\bar{\partial}h_{_j},$ hence 
$$\gr_{_P}(\gk({\frac{\partial}{\partial t}}))=
\sum a_{_{Lj}}\go_{_L}\partial h_{_j},$$
and the theorem is proved.\qed
\brem 
{\em We can also interpret \gr as follows.\\
Fix $\gx\in H^1(T_X),$  then there is a map 
\[
\begin{array}{cccc}
\gt_{_\gx}: & \oplus_kI_k(K_X) & \rightarrow & \oplus_kH^0(K_X^k) \\
& P & \rightarrow & \gr_{_P}(\gx)
\end{array}
\]
Suppose that the canonical map 
$\gi_{_{K_X}}:X\rightarrow\BbP H^0K_X^*=\BbP^{n}$ 
is an embedding, then the polynomials in this $\BbP^{n}$  are 
$\sym H^0(K_X)=\BbC[{\bf x}]$ and the 
ideal of the image $X\simeq\gi_{_{K_X}}(X)\subseteq\BbP^{n}$ is 
$I=\oplus_kI_k(K_X);$ furthermore, if it is projectively normal, 
then its homogeneous coordinate ring 
$S=\frac{\BbC[{\bf x}]}{I}$ coincides with $\oplus_kH^0(K_X^k).$ 
It is easily seen that $\gt_{_\gx}$ is $\BbC[{\bf x}]$--linear, so it factors 
through the quotient, $\gt_{_\gx}:I/I^2\rightarrow S.$\\
Summing up, we have a linear map of graded modules, thus also a natural map 
$$\gt:H^1(T_X)\rightarrow\Hom(I/I^2, S)_0.$$
It is a result of Hilbert scheme theory (see e.g. \cite{s} ch.9) that there 
exists a natural map
$$p_{_0}:\Hom(I/I^2, S)_0\rightarrow H^0(N_{X|\BbP^{n}}),$$ 
where $H^0(N_{X|\BbP^n})$ parametrizes the first--order deformations of $X$ in 
$\BbP^{n}.$ Furthermore, under the hypothesis that $X$ be projectively normal, 
$p_{_0}$ is an isomorphism.\\
By composition, $s:=p_{_0}\circ\gt:H^1(T_X)\rightarrow H^0(N_{X|\BbP^{n}})$
is a section of the normal sequence
$0\rightarrow T_X\rightarrow T_{\BbP^n}\rightarrow N_{X|\BbP^n}\rightarrow0.$}
\erem

\section{Curves}
Let ${\cal X}@>>>B$ be a smooth  family of curves of genus $g\geq3.$ 
The sequence
(\ref{famiglia}) of the previous section has the following natural interpretation in terms of moduli, when considered in degree $k=2.$\\
Let $B$ be any open subset of the moduli space of curves,
outside of the locus of curves with automorphisms. Recall that  over $B$ exists the universal family
$\psi:{\cal C}\rightarrow B$ and the period map 
$\gt:{\cal M}_g\rightarrow{\cal A}_g,$ where ${\cal A}_g$ is the moduli
space of principally polarized abelian varieties of dimension $g,$ becomes an 
embedding when restricted to $B.$\\
In such a situation, there are the following identifications
$$ {\cal I}_2(K_{{\cal C}|B})\simeq
N^*_{\gt(B)|{\cal A}_g},\Hsp
\sym^2(\psi_*K_{{\cal C}|B})\simeq
T^*_{{\cal A}_g}|_{\gt(B)},\Hsp
\psi_*K^2_{{\cal C}|B}\simeq T^*_\gt(B),$$
hence, dualizing (\ref{famiglia}), we obtain the normal sequence
$$0\rightarrow T_{\gt(B)}\rightarrow
T_{{\cal A}_g}|_{\gt(B)}\rightarrow N_{\gt(B)|{\cal A}_g}\rightarrow0.$$
Thus  $\2ff_2:T_{\gt(B)}\rightarrow
{\mathit Hom}({\cal I}_2(K_{{\cal C}|B}),\psi_*K^2_{{\cal C}|B})$
is the same as the 2ff of the sheaf sequence above. Thanks to 
theorem \ref{hg2ff}, we see that \gr is (a factor of) the 2ff of the (local) 
embedding given by the period map for curves.
Karpishpan \cite{k} defines a 2ff of period maps coming from VHS, and asks the 
question, whether, in the case of curves, the 2ff lift the second Wahl map. 
In theorem \ref{wahl} below 
we give a positive answer to this question, by lifting the Wahl map to $\gr.$
It should be remarked that in \cite{cime} p.37-8, a similar lifting of the 
Wahl map is constructed. Indeed, Green states the assertion of 
theorem \ref{wahl} for the case $L=K_C,$ referring for the proof to unpublished
joint work with Griffiths.\\  
For ease of reference, we collect here a few well known facts about Wahl maps 
and Schiffer variations 
(see e.g. \cite{w} and \cite{gr}.)\vspace{\baselineskip}\\
{\em Wahl maps}\vspace{\baselineskip}\\
Let $C$ be a smooth projective curve. Set $S:=C\times C$ and let 
$\gD \subseteq S$ be the diagonal subset. Given a 
line bundle $L$ on $C,$ and $\gl_1,\dots,\gl_r$ a basis of 
$H^0(L)$ as before, define $L_S:=p_1^\ast L\otimes p_2^\ast L,\hsp p_i$ being 
the projection on the i--th factor, $p_i:S\rightarrow C,\hsp i=1,2.$\\
Wahl maps are the natural maps:
$$\gm_{_n}:H^0(S,L_S(-n\gD))\rightarrow 
H^0(S,L_S(-n\gD)|_\gD)\simeq H^0(C,L^2\otimes K_C^n).$$
Especially, $I_2(L)$ can be identified with a subspace of $H^0(S,L_S(-2\gD)).$ 
We are interested in the restriction of $\gm_{_2}$ to $I_2(L),$ i.e. the map 
$$\gm_{_2}:I_2(L)\rightarrow H^0(L^2\otimes K_C^2).$$ 
Also, we recall the local expression of $\gm_{_2}.$ In local coordinates, 
$\gl_i$ is of type
$\gl_i=\phi_i\ell,$ with $\phi_i$ a holomorphic function and $\ell$ a local 
generator of $L;$  then $\sum a_{ij}\gl_i\otimes\gl_j$ is an element of 
$I_2(L)$ iff $\sum a_{ij}\phi_i\phi_j$ is identically zero, and, since the 
$a_{ij}$ are symmetric, also $\sum a_{ij}\dot\gph_i\gph_j=0.$ The local 
expression of $\gm_{_2}$ is 
\begin{equation}
\gm_{_2}\left(\sum a_{ij}\gl_i\otimes\gl_j\right)=
\sum a_{ij}\ddot\gph_i\gph_j\ell^2\otimes dz^2.\label{alloro}
\end{equation}
{\em Schiffer variations}\vspace{\baselineskip}\\
As above, let $L$ be a line bundle over a curve $C,\; \deg L\geq2.$
For any point $P\in C$ consider the exact sequence
$0\rightarrow L^{-1}\rightarrow L^{-1}(P)\rightarrow
L^{-1}(P)|_{_P}\rightarrow0.$\\
The image of the induced natural map 
$\gd:H^0(L^{-1}(P)|_{_P})\rightarrow H^1(L^{-1})$ has dimension one. 
Every generator of $\nim\gd$ is called a Schiffer variation 
of $L$ at $P,$ denoted $\gx_{_P}.$
It is easy to check that, via the Dolbeault isomorphism
$H^1(L^{-1})=H^{0,1}(L^{-1}),\hs \gx_{_P}$ is represented by a form
\begin{equation}
\gth_{_P}=\frac{1}{z}\bar{\partial}b\otimes\ell^\ast,\label{shif}
\end{equation}
where $z$ is a holomorphic coordinate on $C$ around $P,$ $b$ is a bump function
around $P$ and $\ell^*$ is the dual of a local generator of $L.$
\vspace{\baselineskip}\\
{\em A lifting of the Wahl map}
\bt\label{wahl}
The following diagram
\begin{equation}
\begin{CD}
I_2(L) @>\rho>>  \sym^2H^0(L\otimes K_C)\\
@V{\gm_{_2}}VV @VV{\mathbf m}V\\
H^0(L^2\otimes K_C^2) @= H^0(L^2\otimes K_C^2)
\end{CD}\nonumber
\end{equation}
is commutative up to a constant.
\et
The strategy of proof is the following.\\
Given $Q\in I_2(L),$ in order to check $\gm_{_2}(Q)=({\mathbf m}\circ\rho)(Q),$ 
it is enough to evaluate both at every point $P$ in some open subset $U$ of $X.$
$\gm_{_2}(Q)(P)$ is easily computed in terms of (\ref{alloro}). To evaluate 
$({\mathbf m}\circ\rho)(Q)(P)$ we express the dual map ${\mathbf m}^*$ in terms 
of Schiffer variations: namely, if $v_{_P}$ is the evaluation map at  $P,$ then, 
up to a constant, ${\mathbf m}^\ast(v_{_P})=\gx_{_P}\odot\gx_{_P}.$ Thus
$({\mathbf m}\circ\rho)(Q)(P)=(\gx_{_P}\odot\gx_{_P})(\gr(Q)),$ and the right 
hand value is computed making use of the explicit representation (\ref{shif}) 
of the Schiffer variation at $P.$\vspace{\baselineskip}\\  
\pf If $I_2(L)=0,$ there is nothing to prove, so we can suppose $h^0(L)>1$
and $\deg L\geq2.$
For any $P\in C,$ let $v_{_P}$ be the evaluation map at  $P,$ defined on
$H^0(L^2\otimes K_C^2).$ \\
Fix $P_0\in C,$ choose a coordinate $z$ on $C$ and a trivialization of $L,$ 
with local generator $\ell,$ around $P_0.$ Also, let $b$ be a bump function 
around $P_0,$ and let $U\subseteq C$ be an open neighborhood of $P_0$ on which 
$b\equiv1.$ We can suppose that both the  coordinate $z$ on $C$ and the 
trivialization of $L$ are defined on $U.$\\
For all $P\in U,$ via the identification 
$(L^2\otimes K_C^2)_{_P}\simeq\BbC$ coming from the chosen trivialization, 
we can think of $v_{_P}$ as an element of $H^0(L^2\otimes K_C^2)^\ast.$
We want to express its image, under the dual multiplication map 
$${\mathbf m}^\ast:H^0(L^2\otimes K_C^2)^\ast\rightarrow
\sym^2H^0(L\otimes K_C)^\ast,$$
in terms of the Schiffer variation $\gx_{_P}$ of $L$ at $P$ represented by the 
form $\gth_{_P}=\frac{1}{z-z(P)}\bar{\partial}b\otimes\ell^\ast.$
\vspace{1.5\baselineskip}\\{\em Claim}\Hsp 
${\mathbf m}^\ast(v_{_P})=
\frac{1}{(2\pi i)^2}\gx_{_P}\odot\gx_{_P}\in \sym^2H^1(L^{-1}).$
\vspace{\baselineskip}\\
Locally around $P,\hsp\gt\in \sym^2H^0(L\otimes K_C)$ 
has the form 
$\gt=\sum a_{_{ij}}(\phi_{_i}\ell\otimes dz)\otimes (\phi_{_j}\ell\otimes dz).$ 
Thus 
$$v_{_P}({\mathbf m}(\gt))=\sum a_{_{ij}}\phi_{_i}(P)\phi_{_j}(P).$$
On the other hand,  
$$(\gx_{_P}\odot\gx_{_P})(\gt)=
(2\pi i)^2\sum a_{_{ij}}\phi_{_i}(P)\phi_{_j}(P).$$
Now, let $Q=\sum a_{ij}\gl_i\otimes\gl_j\in I_2(L);$ to prove 
$k\cdot\gm_{_2}(Q)=({\mathbf m}\circ\rho)(Q),\hs k$ constant, it is enough to 
show that, for any 
$P\in U,\hsp k\cdot\gm_{_2}(Q)(P)=({\mathbf m}\circ\rho)(Q)(P),$ 
hence, for some constant $h,$ 
$$h\cdot v_{_P}(\gm_{_2}(Q))=(\gx_{_P}\odot\gx_{_P})(\rho(Q)).$$ 
As before write $\gl_i=\phi_i\ell.$ By (\ref{alloro}) 
$$v_{_P}(\gm_{_2}(Q))=\sum a_{_{ij}}\ddot{\phi}_{_i}(P)\phi_{_j}(P)$$
and, by (\ref{last})
\begin{eqnarray}
(\gx_{_P}\odot\gx_{_P})(\rho(Q))&=&\gx_{_P}\cdot\gr_{_Q}(\gx_{_P})\nonumber\\
&=&\left(\frac{1}{z-z(P)}\bar{\partial}b\otimes\ell^\ast\right)\cdot
\sum a_{_{ij}}\gl_{_i}\partial h_{_j}\nonumber\\
&=&\int_C\frac{1}{z-z(P)}\bar{\partial}b
\sum a_{_{ij}}\phi_{_i}\frac{\partial h_{_j}}{\partial z}dz .\nonumber
\end{eqnarray}
\vspace{\baselineskip}\\
Write $\Psi(z):=\sum a_{_{ij}}\phi_{_i}\frac{\partial h_{_j}}{\partial z},$ 
then the value of the last integral is $2\pi i\Psi(P).$\\
To evaluate $\Psi(P)$ we proceed as follows: 
$\gth_{_P}\gl_{_i}=\frac{\phi_{_i}}{z-z(P)}\bar{\partial}b,$
so, in $C-\{P\},$ we have the equality (cf.(\ref{decomposition}))
$\gga_{_i}+\bar{\partial}h_{_i}=
\bar{\partial}\left(\frac{b\phi_{_i}}{z-z(P)}\right),$
hence $\gga_{_i}=\bar{\partial}g_{_i},$ with $g_{_i}=\frac{b\phi_{_i}}{z-z(P)}-h_{_i}.$\\ 
Define $\eta_{_i}:=\partial g_{_i}.$
\bl\label{lemmino}{\em (cf. \cite{p} 4.8)}
The $\eta_{_i}$ are all proportional, hence
$$\gr_{_Q}(\gx_{_P})=\eta\sum a_{_{ij}}b_{_i}\gl_{_j}\in H^0(L\otimes K_C),$$
where \get is a differential of second kind, multiple of the $\eta_{_i},$ 
having only a double pole at $P,$
and $b_{_i}=\frac{\get_{_i}}{\get}(P)$ are constants.
\el
\pf {\em Step 1}\Hsp 
$\sum a_{_{ij}}\gl_{_i}\partial h_{_j}=-\sum a_{_{ij}}\gl_{_i}\eta_{_j}.$
\\
In the first place, note that $\eta_{_i}$ is holomorphic in  $C-\{P\}.$
Indeed, $\gga_{_i}$ is harmonic, thus  
$$\bar{\partial}\eta_{_i}=\bar{\partial}\partial g_{_i}=
-\partial\bar{\partial}g_{_i}=-\partial\gga_{_i}=0.$$ 
Hence,
$\sum a_{_{ij}}\gl_{_i}\partial(\frac{b\phi_{_j}}{z-z(P)})=\sum a_{_{ij}}\gl_{_i}\partial h_{_j}+
\sum a_{_{ij}}\gl_{_i}\eta_{_j}$  is holomorphic on $C-\{P\}$ and, in a neighborhood of 
$P$ where $b\equiv1,$ it has the form 
$$\left(\sum a_{_{ij}}\phi_{_i}\frac{\partial}{\partial z}
\left(\frac{\phi_{_j}}{z-z(P)}\right)\right)\ell\otimes dz=
\sum a_{_{ij}}\phi_{_i}\left(-\frac{\phi_{_j}}{(z-z(P))^2}+
\frac{\dot{\phi_{_j}}}{z-z(P)}\right)\ell\otimes dz=0,$$
because $\sum a_{_{ij}}\phi_{_i}\phi_{_j}=\sum a_{_{ij}}\phi_{_i}\dot{\phi_{_j}}=0,$ so 
$\sum a_{_{ij}}\gl_{_i}\partial(\frac{b\phi_{_j}}{z-z(P)})$ 
is identically zero. 
\vspace{1.5\baselineskip}\\{\em Step 2}\Hsp 
$\eta_{_i}\in H^0(K_C(2P)).$ 
\\
On $U,$ where $b\equiv1, \hs\eta_{_i}$ has the form 
\begin{equation}
\eta_{_i}=\left(-\frac{\phi_{_i}(P)}{(z-z(P))^2}+f_{_i}(z)\right)dz\label{eta}
\end{equation}
with $f_{_i}(z)$ a holomorphic function, hence $\eta_{_i}$ is a meromorphic
form, with a double pole at $P.$
\vspace{1.5\baselineskip}\\{\em Step 3}\Hsp 
$\eta_{_i}$ are proportional. 
\\
If $C$ has genus $g=0,$ then $h^0(K_C(2P))=1.$\\
If $g\geq1,$ by definition, 
$\eta_{_i}+\gga_{_i}=\partial g_{_i}+\bar{\partial}g_{_i}=dg_{_i},$ so 
$[\eta_{_i}]=-[\gga_{_i}]\in H^1(C-\{ P\},\BbC),$
hence $\eta_{_i}\in H^0(K_C(2P))\cap H^{0,1}(C),$ via the inclusion 
$H^0(K_C(2P))\hookrightarrow H^1(C-\{ P\},\BbC)\simeq H^1(C,\BbC).$ \\
Now, $\dim H^0(K_C(2P))\cap H^{0,1}(C)=1:$ it is a consequence of 
$h^0(K_C(2P))=g+1$ and $H^0(K_C)\subseteq H^0(K_C(2P)).$
It follows that the $\eta_{_i}$ are all proportional.\qed  
To finish the proof of the theorem, take $\eta$ 
to be the only differential
having local expression $\eta=(-\frac{1}{(z-z(P))^2}+f(z))dz,\; f(z)$
holomorphic. So $\eta_{_i}=\phi_{_i}(P)\eta,$ and 
we see that $\Psi(z),$ given by the local expression of 
$-\sum a_{_{ij}}\gl_{_i}\eta_{_j},$ has the form 
$\Psi(z)=
\frac{1}{2}\sum a_{_{ij}}\ddot{\phi_{_i}}(P)\phi_{_j}(P)+\circ(z-z(P)).$
Thus, 
$$\Psi(P)=\frac{1}{2}\sum a_{_{ij}}\ddot{\phi_{_i}}(P)\phi_{_j}(P),$$ 
hence $v_{_P}(\gm_{_2}(Q))=(\gx_{_P}\odot\gx_{_P})(\rho(Q)),$ 
up to a constant factor.\qed
\brem\label{lastalimone} 
{\em Recall that, if $S$ is a subbundle of the hermitian bundle $E,$ and 
$\nabla_E$ and $\nabla_S$ are the metric connections, then
the 2ff of the embedding $S\hookrightarrow E$ gives information on the 
curvature of $S,$ because of the relation
$$\nabla_S=\nabla_E-\2ff.$$
So, the lemma above makes possible explicit computations about the curvature 
of the moduli space of curves.}
\erem
\section{Generalizations}
{\em Pairs of vector bundles}\vspace{\baselineskip}\\
Let $E$ and $F$ be vector bundles on a smooth projective variety $X.$ As in the
previous section, define $Y:=X\times X,$ with $p_i:Y\rightarrow X,\hsp i=1,2,$
projections on the i--th factor and $\gD \subseteq Y$ the diagonal subset.
Tensoring the exact sequence
$0\rightarrow {\cal I}_\gD^{n+1}\rightarrow{\cal I}_\gD^n\rightarrow
{\cal I}_\gD^{n}/{\cal I}_\gD^{n+1}\rightarrow 0$ 
by $E\Boxtimes F:=p_1^\ast E\otimes p_2^\ast F,$ and taking cohomology, we get
$$0\rightarrow H^0(Y, E\Boxtimes F\otimes{\cal I}_\gD^{n+1})\rightarrow  
H^0(Y, E\Boxtimes F\otimes{\cal I}_\gD^n)\rightarrow  
H^0(X, E\otimes F\otimes\sym^n \gO^1_X).$$
Following the notation of \cite{pa}, we define the {\em k-th module of relations
of E and F} as $R_k(E,F):=H^0(Y, E\Boxtimes F\otimes{\cal I}_\gD^k)
\subseteq H^0(Y, E\Boxtimes F)\simeq H^0(X,E)\otimes H^0(X,F).$\\
Note also that, when $E=F, \; I_2(E)$ is a 
submodule of $R_2(E,E).$\\ 
We now extend the definition of \gr to $R_2(E,F),$ to obtain a map, still 
denoted by \gr,
$$\gr: R_2(E,F)\rightarrow
\text{\normalshape Hom}(H^{p,q}(X, E^*), H^{p+1,q-1}(X, F)).$$
Let $\gl_{_i}, i=1,\ldots,s,$ and $\gm_{_j}, j=1,\ldots,t,$ be bases of 
$H^0(X, E)$ and $H^0(X, F)$ respectively. Because of the inclusion 
$R_2(E,F)\subseteq H^0(X, E)\otimes H^0(X, F),$ an element $P\in R_2(E,F)$
can be written as $P=\sum a_{_{ij}}\gl_{_i}\otimes \gm_{_j}.$ If 
$\gx\in H^{p,q}(E^*)$ and  $\gth\in A^{p,q}(E^*)$ is a Dolbeault 
representative of \gx, then $\gth\gl_{_i}\in A^{p,q}(X)$ are 
$\bar{\partial}$--closed forms having harmonic decompositons
$\gth\gl_{_i}=\gga_{_i}+\bar{\partial}h_{_i}.$ We define $\gr_{_P}(\gx)$ 
as the Dolbeault cohomology class of the $(p+1,q-1)$--form
$$\gs_{_P}(\gth):=\sum_{_{ij}}a_{_{ij}}\partial h_{_i}\otimes\gm_{_j}.$$
\bt
The map 
$$\begin{array}{cccc}
\gr: & R_2(E,F) & \rightarrow & 
\text{\normalshape Hom}(H^{p,q}(E^*), H^{p+1,q-1}(F))\\
& P & \rightarrow & (\gx\mapsto\gr_{_P}(\gx))
\end{array}
$$
is well defined and linear.
\et
\pf The proof runs along the same lines of that of \ref{uno}.\\
(i)\Hsp $\gs_{_P}(\gth)$ is $\bar{\partial}$--closed.\\
Clearly, 
$\bar{\partial}\gs_{_P}(\gth)=
\bar{\partial}\left(\sum a_{_{ij}}\partial h_{_i}\otimes\gm_{_j}\right)=
-\sum a_{_{ij}}\partial(\gth\gl_{_i})\otimes\gm_{_j}$ 
as before, so we perform again a local computation.\\
Let $U$ be an open subset of $X$ on which $E$ and $F$ are both locally trivial,
and let $\ell_{_i},i=1,\ldots,q,$ and $m_{_j},j=1,\ldots,r,$ be local basis of 
$E$ and $F$ respectively, with $\ell_{_i}^*,i=1,\ldots,q,$ the dual basis of
$E^*$, then locally $\gl_{_i}=\sum_{_l}\gga_{_{il}}\ell_{_l}$ and 
$\gm_{_j}=\sum_{_k}\gb_{_{jk}}m_{_k},$ with $\gga_{_{il}}, \gb_{_{jk}}$
holomorphic functions. An element 
$P=\sum a_{_{ij}}\gl_{_i}\otimes \gm_{_j}\in R_2(E,F),$ being a section of a 
bundle twisted by ${\cal I}_\gD^2,$ vanishes on $X$ to the second order, which
locally translates into $\sum_{_{ij}}a_{_{ij}}\gga_{_{il}}\gb_{_{jk}}=0$ and 
$\sum_{_{ij}}a_{_{ij}}\gb_{_{jk}}\partial\gga_{_{il}}=0$ for all $l, k.$
Also, $\gth\in A^{p,q}(E^*)$ on $U$ is of the form 
$\gth=\sum_{_l}\go_{_l}\ell_{_l}^*,$ with $\go_{_l}\in A^{p,q}(X),$ hence 
$\gth\gl_{_i}=\sum_{_l}\gga_{_{il}}\go_{_l}.$ It follows that locally
\begin{eqnarray}
\sum_{_{ij}} a_{_{ij}}\partial(\gth\gl_{_i})\otimes\gm_{_j}&=& 
\sum_{_{ij}} a_{_{ij}}\partial\left(\sum_{_l}\gga_{_{il}}\go_{_l}\right)\otimes
\sum_{_k}\gb_{_{jk}} m_{_k}\nonumber\\
&=&\sum_{_{ijkl}} a_{_{ij}}\gb_{_{jk}}(\partial\gga_{_{il}}\wedge\go_{_l}+
\gga_{_{il}}\wedge\partial\go_{_l})\otimes m_{_k}\nonumber\\
&=&\sum_{_{kl}}
\left(\sum_{_{ij}}a_{_{ij}}\gb_{_{jk}}\partial\gga_{_{il}}\right)\wedge\go_{_l}
\otimes m_{_k}+\nonumber\\
&&\sum_{_{kl}}\left(\sum_{_{ij}}a_{_{ij}}\gga_{_{il}}\gb_{_{jk}}\right)
\wedge\partial\go_{_l}\otimes m_{_k}\nonumber\\
&=&0\nonumber
\end{eqnarray}
thus $\gs_{_P}(\gth)$ is $\bar{\partial}$--closed.\\
(ii)\Hsp $\gr_{_P}(\gx)$ does not depend on the choice of $\gth.$\\ 
The argument is completely similar to that of \ref{uno}.\qed
\brem
{\em By Serre duality, the range of the map \gr is 
$\text{\normalshape Hom}(H^{p,q}(E^*), H^{p+1,q-1}(F))=
(H^{p,q}(E^*))^*\otimes H^{p+1,q-1}(F)=H^{n-p,n-q}(E)\otimes H^{p+1,q-1}(F).$
Also, it is a consequence of Kunneth formula that
$H^{r,s}(X\times X, E\Boxtimes F)=
\oplus_{\stackrel{i+j=r}{\scriptscriptstyle h+k=s}} 
H^{i,h}(X,E)\otimes H^{j,k}(X,F).$\\
Thus, adding together all the maps \gr, we have a 
natural map (still denoted by \gr) 
$$\gr:H^0(X\times X, E\Boxtimes F\otimes{\cal I}_\gD^2)\rightarrow
H^{n-1}(X\times X, E\Boxtimes F\otimes\gO^{n+1}_{X\times X}).$$}
\erem
{\em Harmonic bundles}
\vspace{\baselineskip}\\
We collect here some definitions and known facts about Higgs and harmonic
bundles (cf. \cite{si}.)\\
Let $X$ be a compact K\"{a}hler manifold of dimension $n,$ with K\"{a}hler
form $\go.$ A Higgs field (or Higgs bundle) is a pair $(E,\phi),$ with $E$ 
a holomorphic vector
bundle and $\phi:E\rightarrow E\otimes\gO_X^1$ a holomorphic map such that 
$\phi\wedge\phi=0.$ Associated to \gph there is the operator 
$D'':=\bar{\partial}+\phi:{\cal A}^0(E)\rightarrow{\cal A}^1(E),$ with $$D''(f\gs)=\bar{\partial}f\cdot\gs+fD''\gs,\Hsp {D''}^2=0.$$
The Dolbeault cohomology $H^*_{{\mathit Dolb}}(E)$ is defined as the 
hypercohomology of the complex
$$E\otimes\gO_X^{\bullet}\Hsp 
E@>\phi>>E\otimes\gO_X^1@>\phi>>E\otimes\gO_X^2@>\phi>>\ldots$$
$H^*_{{\mathit Dolb}}(E)$ is isomorphic to the cohomology of the complex 
${\cal A}^i(E)@>{D''}>>{\cal A}^{i+1}(E).$ Note that 
$D''$ defines a different holomorphic structure on $E.$ 
If $E$ is endowed with a hermitian metric $H,$ define $D'_H$ as the operator 
for which $D=D'_H+D''$ is the hermitian 
connection, with respect to the holomorphic structure of $E$ associated to  
$D''.$ When $D$ is flat, $E$ is called a harmonic bundle 
(and $H$ is a harmonic metric.)\\
A fundamental result in the theory of Higgs and harmonic bundles is the 
following.
\bt\label{higgs}
{\em (cf. \cite{si} theorem 1)}\\
A Higgs bundle has a harmonic metric if and only if it is polystable
(i.e. direct sum of stable Higgs bundles having the same slope) and
$c_{_1}(E)[\go]^{n-1}=c_{_2}(E)[\go]^{n-2}=0.$\\
Conversely, a flat bundle (with a metric) comes from a Higgs bundle if and only if it is semisimple.
\et
The following hold for harmonic bundles:\\
(i) The  K\"{a}hler identities.\\
(ii) The associated harmonic decomposition
$$\begin{array}{rcl}
{\cal A}^{\bullet}(E)&=&{\cal H}^{\bullet}(E)\oplus\nim D\oplus\nim D^*\\
&=&{\cal H}^{\bullet}(E)\oplus\nim D''\oplus\nim {D''}^*,
\end{array}$$
${\cal H}(E)$ being the kernel of the laplacian operator
$\gD=DD^*+D^*D=2(D''{D''}^*+{D''}^*D'').$\\
(iii) The principle of two types
$$\nker(D'_H)\cap\nker(D'')\cap(\nim(D'_H)+\nim(D''))=\nim(D'_HD'').$$
These properties are all one needs to generalize the construction of the map $\gr.$
Let $(E,H)$ be a harmonic bundle, with $D''=\bar{\partial}+\phi.$ Then, for any 
line bundle $L,$ $D''$ defines structures of Higgs bundles on both 
$E\otimes L$ and $E\otimes L^{-1},$ with associated cohomology
$H^*_{{\mathit Dolb}}(E\otimes L)$ and $H^*_{{\mathit Dolb}}(E\otimes L^{-1}).$ 
\bt
Let $(E,H)$ be a harmonic vector bundle on a compact
K\"{a}hler manifold $X,$ and let $L$ be any line bundle on $X.$ 
Then is well--defined the map
$$\begin{array}{rcl}
I_2(L)\otimes H^k_{{\mathit Dolb}}(E\otimes L^{-1})&\rightarrow&
H^k_{{\mathit Dolb}}(E\otimes L)\\
(Q, \ga)&\rightarrow&[\sum a_{_{ij}}\gl_{_i}D'_Hh_{_j}]
\end{array}$$
where:\\
(i)\hsp $\{ \gl_{_i}\} $ is any basis of $H^0(L),$ so that  $Q$ can be written
as $Q=\sum a_{_{ij}}\gl_{_i}\otimes\gl_{_j},$ and\\
(ii)\hsp $h_{_j}$ is given by the harmonic decomposition $\gl_{_j}\tilde{\ga}=\gga_{_j}+D''h_{_j},$\hsp $\tilde{\ga}$ being a form 
representing \ga, i.e. $\ga=[\tilde{\ga}],$ with
$\tilde{\ga}\in{\cal A}^i(E\otimes L^{-1})$ and $D''\tilde{\ga}=0.$
\et
\pf The proof is completely analogous to that of 
proposition--definition \ref{uno}.\qed
The simplest case is that of a polystable vector bundle with
$c_{_1}(E)[\go]^{n-1}=c_{_2}(E)[\go]^{n-2}=0$---in other terms,i.e. $D''=\bar{\partial}$ and 
$H^k_{{\mathit Dolb}}(E)=\oplus_{p=0}^kH^{k-p}(E\otimes\gO_X^p),$ 
also \gr is a map
$$I_2(L)\rightarrow\oplus_p H^{n-k+p}(E^*\otimes L\otimes\gO_X^{n-p})
\otimes H^{k-p-1}(E\otimes L\otimes\gO_X^{p+1}).$$
Especially, for any degree zero line bundle $M$ on a smooth curve $C,$ there
exists a harmonic metric $H$ on $M,$ with metric connection $D_H$
that decomposes as $D_H=D_H'+\bar{\partial}.$ Thus we have the map 
$$I_2(L)@>\gr>>
H^0(M^{-1}\otimes L\otimes K_C)\otimes H^0(M\otimes L\otimes K_C).$$ 
By means of the multiplication map
$$H^0(M^{-1}\otimes L\otimes K_C)\otimes H^0(M\otimes L\otimes K_C)
@>m>>H^0(L^2\otimes K_C^2),$$
we have the following generalization of theorem \ref{wahl}.
\bt
The diagram
\begin{equation}
\begin{CD}
I_2(L) @>\rho>>
H^0(M^{-1}\otimes L\otimes K_C)\otimes H^0(M\otimes L\otimes K_C)\\
@V{\gm_{_2}}VV @VVmV\\
H^0(L^2\otimes K_C^2) @= H^0(L^2\otimes K_C^2)
\end{CD}\nonumber
\end{equation}
is commutative up to a constant.
\et
\pf The proof, that uses the operators $D,\: D_H'$ and $\bar{\partial}$ in the r\^oles of $d,\: \partial$ and $\bar{\partial}$ respectively, is analogous to that of theorem \ref{wahl} but for the details noted below.\\
(i)\Hsp We suppose that on $U$ there exists also a trivialization of $M,$ with 
local generator $\nu.$ Then, in the claim, the Schiffer variation 
$\gx_{_P}\in H^1(M\otimes L^{-1})$ is the one represented by the form 
$\gth_{_P}=\frac{1}{z-z(P)}\bar{\partial}b\otimes\gn\otimes\ell^\ast.$\\
(ii)\Hsp The metric is represented on $U$ by a scalar function, still denoted by
$H.$ Hence $D_H':{\cal A}^0(M)\rightarrow{\cal A}^{0,1}(M)$ locally is 
$$D_H'(f\nu)=\left(\partial f+f\frac{\partial H}{H}\right)\otimes\gn.$$
Writing $h_{_i}=l_{_i}\gn,$ then $\gr_{_Q}(\gx)$ is represented by the form 
$\Psi(z)\gn\otimes\ell\otimes dz,$ with 
$$\Psi(z)=\sum a_{_{ij}}\phi_{_i}
\left(\frac{\partial l_{_j}}{\partial z}+
l_{_j}\frac{\partial\log H}{\partial z}\right).$$
(iii)\Hsp Steps 1 and 2 of lemma \ref{lemmino} carry through the present 
situation, with the local expression (\ref{eta}) for the form $\get_{_i}$ now
becoming
\begin{equation}
\get_{_i}=\left(-\frac{\phi_{_i}(P)}{(z-z(P))^2}+
\frac{\partial\log H}{\partial z}(z(P))
\frac{\phi_{_i}(P)}{z-z(P)}+f_{_i}(z)\right)\nu\otimes dz.\label{nuova}
\end{equation}
To prove step 3, we argue as follows.\\
Assume that $P$ is not a base point for $L,$ i.e. not all $\phi_{_i}(P)=0.$
As a consequence of (\ref{nuova}), we have that 
$\phi_{_j}(P)\eta_{_i}-\phi_{_i}(P)\eta_{_j}$ is a $D_H'$--harmonic form defined on all $C.$
Now $\eta_{_i}+\gga_{_i}=D_H'g_{_i}+\bar{\partial}g_{_i}=D_Hg_{_i},$ so 
$$\phi_{_j}(P)\eta_{_i}-\phi_{_i}(P)\eta_{_j}+\phi_{_j}(P)\gga_{_i}-\phi_{_i}(P)\gga_{_j}=
D_H\left(\phi_{_j}(P)g_{_i}-\phi_{_i}(P)g_{_j}\right).$$
The equality above shows that the $D_H$--harmonic form on the left hand side is $D_H$--exact, hence it is zero, because of the harmonic decomposition.
Especially, its $D_H'$-part is zero, hence
$\phi_{_j}(P)\eta_{_i}=\phi_{_i}(P)\eta_{_j}.$ So the $\get_{_i}$'s are 
proportional and the final computation of the proof can still be performed.\qed


\section*{}
\begin{tabbing}
e--mail addresses: \= 
{\tt colombo{\char`@}mat.unimi.it}\\
\> {\tt pirola{\char`@}dimat.unipv.it}\\
\> {\tt tortora{\char`@}mat.unimi.it}\\
\end{tabbing}



\begin{thebibliography}{m3n7}

\bibitem[CGGH]{ivhs} J. Carlson, M. Green, P. Griffiths, J. Harris, 
{\em Infinitesimal variations of Hodge structure (I),} 
Compositio Math. {\bf 50} (1983), 109-205.

\bibitem[DGMS]{dgms} P. Deligne, P. Griffiths, J. Morgan, D. Sullivan, 
{\em Real homotopy theory on K\"ahler manifolds,} 
Invent. Math. {\bf 29} (1975), 245-274.

\bibitem[G]{cime} M. Green, {\em Infinitesimal methods in Hodge theory,}
Algebraic cycles and Hodge theory, Lecture Notes in Math., vol. 1594 (1994), 
1-92.

\bibitem[Gr]{gr} P. Griffiths,
{\em Infinitesimal variations of Hodge structure (III),} 
Compositio Math. {\bf 50} (1983), 267-324.

\bibitem[GH1]{pag} P. Griffiths, J. Harris, 
{\em Principles of algebraic geometry,} New York, John Wiley, 1978.

\bibitem[GH2]{gh} P. Griffiths, J. Harris, {\em Algebraic geometry and local 
differential geometry,} 
Ann. Sci. \`{E}cole Norm. Sup. (4) {\bf 12} (1979), 355-432.

\bibitem[K]{k} Y. Karpishpan, {\em On higher--order differentials of the period 
map,} Duke Math. J. {\bf 72} (1993), 541-571.

\bibitem[L]{l} J. Landsberg, {\em On second fundamental form of projective
varieties,} Invent. Math. {\bf 117} (1994), 303-315.

\bibitem[P]{pa} R. Paoletti, {\em Generalized Wahl maps and adjoint line 
bundles on a general curve,} Pacific J. Math. {\bf 168} (1995), 313-334.

\bibitem[Pi]{p} G.P. Pirola, {\em The infinitesimal variation of the spin 
abelian differentials and periodic minimal surfaces,}  
Comm. Anal. Geom. {\bf 6} (1998), 393-426.

\bibitem[S]{s} E. Sernesi, {\em Topics on families of projective schemes,}
Queen's Papers no.73, Kingston, Ontario, 1986.

\bibitem[Si]{si} C. Simpson, {\em Higgs bundles and local systems,}
Publ. Math. I.H.E.S. {\bf 75} (1992), 5-95.

\bibitem[W]{w} J. Wahl, {\em Gaussian maps on algebraic curves,} 
J. Diff. Geom. {\bf 32} (1990), 77-98.

\end{thebibliography}
\end{document}